\documentclass[11pt]{amsart}
\usepackage{amssymb,amsmath,amsthm}
\usepackage[margin=1in]{geometry}
\usepackage[hidelinks]{hyperref}

\theoremstyle{plain}
\newtheorem{theorem}{Theorem}[section]
\newtheorem{proposition}[theorem]{Proposition}
\newtheorem{lemma}[theorem]{Lemma}
\newtheorem{corollary}[theorem]{Corollary}
\newtheorem*{theoremA}{Theorem A}
\newtheorem*{theoremB}{Theorem B}
\newtheorem*{theoremC}{Theorem C}
\theoremstyle{definition}

\theoremstyle{remark}
\newtheorem{remark}[theorem]{Remark}

\newcommand{\eps}{\varepsilon}
\newcommand{\N}{\mathbb{N}}
\newcommand{\Z}{\mathbb{Z}}
\newcommand{\Q}{\mathbb{Q}}

\newcommand{\Qbar}{\overline{\mathbb{Q}}}
\newcommand{\pT}{p_T}
\newcommand{\cC}{\mathcal{C}}
\newcommand{\lean}[1]{\texttt{#1}}

\makeatletter
\def\@settitle{\begin{center}%
  \baselineskip14\p@\relax
  \bfseries\Large
  \@title
  \end{center}%
}
\def\@setauthors{%
  \begingroup
  \def\thanks{\protect\thanks@warning}%
  \trivlist
  \centering
  \normalsize
  \@topsep30\p@\relax
  \advance\@topsep by -\baselineskip
  \item\relax
  \andify\authors
  \authors
  \endtrivlist
  \endgroup
}
\makeatother

\title[The interval result for a set of weights]{The Evertse--Ferretti interval result for a set
of weights, and the exponent $1/2$ in three complexity bounds}
\author{Ralf Stephan}
\thanks{Institute for Globally Distributed Open Research and Education (IGDORE) \\
MSC 11J68, 11J87, 11A63, 68R15}
\date{October 2026}

\begin{document}

\begin{abstract}
The interval result of Evertse and Ferretti (2013) concerns a tuple $L$ of linear forms and one
tuple $c$ of weights. It says that the parameters $Q$ for which the twisted height
$H_{L,c,Q}$ has a small value outside an exceptional subspace $T(L,c)$ lie in
$O(\delta^{-2}\log\delta^{-1})$ intervals $[Q_i,Q_i^{\omega})$. We prove this for a whole set
$\cC$ of weights with the same forms. At most $J$ subspaces $T(L,c)$ occur, where $J$ depends
only on the dimension and on the number of places, and one family of
$O(\delta^{-2}\log\delta^{-1})$ intervals serves every $c\in\cC$ at once. Neither count depends
on $\cC$. In
applications of Bugeaud--Evertse type the solutions are sorted into about $\eps^{-1}$ classes
with the same forms and different weights, and one application now serves all classes. We
deduce three complexity bounds. For a real algebraic irrational number $\xi$, the number
$p(n,\xi,b)$ of blocks of length $n$ in its base-$b$ expansion satisfies
$\limsup p(n,\xi,b)/(n(\log n)^\eta)=\infty$ for every $\eta<1/2$; Bugeaud and Evertse proved
this for $\eta<1/11$. The same holds for the continued fraction expansion of an algebraic number
of degree at least three, and for the steering word of $(3/2)^n$. We explain why $1/2$ is the
limit of the method. All results are machine-verified in Lean~4.
\end{abstract}

\maketitle

\section{Introduction}

\subsection*{The interval result}
Let $K$ be a number field, $L=(L_i^{(v)})$ a tuple of linear forms in $n$ variables, one
$n$-tuple of linearly independent forms for each place $v$ of $K$, and $c=(c_{iv})$ a tuple of
real weights. The twisted height $H_{L,c,Q}(x)=\prod_v\max_i\|L_i^{(v)}(x)\|_vQ^{-c_{iv}}$ of
Evertse and Ferretti \cite[(2.13)]{EF13} measures how well a point $x\in\Qbar^n$ approximates
the system of inequalities $\|L_i^{(v)}(x)\|_v\le Q^{c_{iv}}$ at the level $Q$. Their interval
result \cite[Thm.~2.3]{EF13} says that there is a proper subspace $T=T(L,c)$ of $\Qbar^n$ and
\[
  m_0\asymp\delta^{-2}\log\delta^{-1}
\]
reals $Q_1<\dots<Q_{m_0}$ such that every $Q$ with a point $x\notin T$ of twisted height at most
$\Delta_L^{1/n}Q^{-\delta}$ lies in $[1,C_0)$ or in one of the intervals $[Q_h,Q_h^{\omega_0})$,
$\omega_0=\delta^{-1}\log3R$. Theorem~3.3 of \cite{EF13} is the corresponding statement for
systems of inequalities. It improves the interval result of Faltings and W\"ustholz \cite{FW94}.

The theorem concerns one tuple $c$. In applications the same forms often come with many weights.
This happens in Bugeaud and Evertse's method for complexity bounds \cite{BE08}: the solutions
are sorted into about $\eps^{-1}$ classes, and the classes have the same forms but different
weights. With one application per class, the classes cost a factor $\eps^{-1}$.

\subsection*{The interval result for a set of weights}
Our main result is Theorem~\ref{thm:F} below. Let the forms $L$ be fixed and let $\cC$ be a set
of weights satisfying EF13's normalisations. Assume that the weights are bounded by a common
$\gamma_v$ at each place, with $\sum_v\gamma_v\le1$, and that they all vanish outside a fixed
set $S_1$ of places. Then:
\begin{itemize}
\item the exceptional subspaces $T(L,c)$, $c\in\cC$, are at most
  $J=(n+1)\,2^{2^n\#S_1}$ subspaces;
\item one family of at most $J\,n^{n\#S_1}m(n,R,\delta)$ intervals $[Q_h,Q_h^{\omega_0})$ serves
  every $c\in\cC$, where $m(n,R,\delta)\ll_{n,R}\delta^{-2}(1+\log\delta^{-1})$ is a number of
  the same kind as $m_0$.
\end{itemize}
Neither count depends on $\cC$, which may be infinite. The two ingredients are these.
\begin{itemize}
\item Evertse and Ferretti's version of Roth's method treats the auxiliary polynomial block by
  block. So the weights may change from one block to the next, as long as they share the bound
  $\gamma_v$. Evertse and Ferretti point out that their approach ``allows us to handle exponents
  $\hat c_{iw}(Q)$ that vary with $Q$'' \cite[p.~549]{EF13}. There this concerns auxiliary
  exponents at one place in the proof of the semistable case; we use the same feature for a
  whole set of weights.
\item The subspace $T(L,c)$ is determined by a finite combinatorial datum, its \emph{weight
  type}, and the uniqueness of $T(L,c)$ turns this into the explicit bound $J$. That
  $T(L,c)$ ranges over a finite collection independent of $c$ is known (see below); what is new
  is that the intervals can be chosen uniformly.
\end{itemize}
Section~\ref{sec:systems} carries Theorem~\ref{thm:F} over to systems of inequalities
(Theorem~3.3 of \cite{EF13} for a set of weights). For solutions of geometrically growing height
it gives $O(\delta^{-2}\log^2\delta^{-1})$ solutions outside a bounded number of subspaces,
uniformly over the whole set (Corollary~\ref{cor:doubling}).

\subsection*{Complexity bounds}
Let $b\ge2$ be an integer and $\xi$ a real number with $b$-ary expansion
$\xi=\lfloor\xi\rfloor+\sum_{j\ge1}a_jb^{-j}$. The block complexity
$p(n,\xi,b)=\operatorname{Card}\{a_{k+1}\cdots a_{k+n}:k\ge0\}$ \cite{BE08} is $b^n$ for every $n$
when $\xi$ is normal to base $b$, which is expected for every algebraic irrational $\xi$. The
known lower bounds are barely above linear. Adamczewski and Bugeaud \cite{AB07} proved
$p(n,\xi,b)/n\to\infty$. Bugeaud and Evertse \cite[Thm.~2.1]{BE08} proved
\begin{equation}\label{eq:BE}
  \limsup_{n\to\infty}\frac{p(n,\xi,b)}{n(\log n)^\eta}=+\infty\qquad\text{for every }\eta<1/11.
\end{equation}
In \cite{Ste26b} we raised $1/11$ to $1/4$, together with the analogous statements for continued
fractions, announced by Bugeaud \cite[\S6]{Bug13}, and for the steering word of $(3/2)^n$
\cite[Problem~5.2]{St26}. The exponent $1/4$ is $1/(1+3)$: $\eps^{-1}$ classes times the
$\delta^{-3}$ subspaces of Evertse and Ferretti's count \cite[Thm.~3.1]{EF13}. Here we prove the
following.

\begin{theoremA}
Let $b\ge2$ be an integer and $\xi$ a real algebraic irrational number. For every real $\eta<1/2$
and every real $C$,
\[
  p(n,\xi,b)>C\,n(\log n)^\eta\qquad\text{for infinitely many }n.
\]
\end{theoremA}

\begin{theoremB}
Let $a=(a_\ell)_{\ell\ge1}$ be a sequence of positive integers such that
$\alpha=[0;a_1,a_2,\dots]$ is algebraic of degree at least three, and let
$p(n,a)=\operatorname{Card}\{a_{\ell+1}\cdots a_{\ell+n}:\ell\ge0\}$. For every $\delta$ with
$0\le\delta<1/2$ and every real $C$,
\[
  p(n,a)>C\,n(\log n)^\delta\qquad\text{for infinitely many }n.
\]
\end{theoremB}

\begin{theoremC}
Write $(3/2)^n=m_n+\eps_n$ with $m_n$ the nearest integer, and let $\pT(k)$ be the subword
complexity of the word $T=(2m_{n+1}-3m_n)_{n\ge0}$. For every real $\eta<1/2$ and every real $C$,
\[
  \pT(k)>C\,k(\log k)^\eta\qquad\text{for infinitely many }k.
\]
\end{theoremC}

Each statement is equivalent to $\limsup p(n)/(n(\log n)^\eta)=+\infty$. In Theorem~B, $p(n,a)$
may be infinite; this happens exactly when $a$ is unbounded, and then there is nothing to prove.

The three proofs follow \cite[\S\S6--7]{BE08} and \cite{Ste26b}. A complexity bound
$p(n)\le cn(\log n)^u$ gives $N$ integer points of doubly exponentially growing height, each
satisfying a system of inequalities with exponent $\eps=N^{-v+o(1)}$, $u<v$. Two changes give
$1/2$.
\begin{itemize}
\item We use the interval form instead of the count of subspaces. An interval
  $[Q_h,Q_h^{\omega_0})$ contains only $O(\log\omega_0)$ points of doubly exponentially growing
  height, so the $\delta^{-2}$ intervals hold $O(\delta^{-2}\log^2\delta^{-1})$ points. With
  one application per class this already gives $\eps^{-1}\cdot\eps^{-2}$, i.e.\ every exponent
  below $1/3$.
\item With Theorem~\ref{thm:F} one application serves all classes, and the factor $\eps^{-1}$
  of the classes disappears. The subspaces are now a bounded number, each holding
  $O(\log\eps^{-1})$ points by the arguments of \cite{Ste26b}.
\end{itemize}
Hence $N\le\eps^{-2+o(1)}=N^{2v+o(1)}$, which is impossible for large $N$ when $2v<1$.

\subsection*{The limit of the method}
The exponent $1/2$ is the end of this road. The remaining $\delta^{-2}$ is the number of blocks of
the auxiliary polynomial in Roth's method. It comes from the probabilistic estimate (Hoeffding's
inequality) that controls the index. The quantitative Roth theorem has the same $\delta^{-2}$
\cite[Prop.~1.2]{Ev10}. A sequence of doubly exponentially growing height can place one point
into each of the $m_0$ intervals, so no bound of this shape gives fewer than $\delta^{-2}$
points (Remark~\ref{rem:limit}).

\subsection*{Prior work}
That the exceptional subspaces come from a finite collection depending only on the forms, and
not on the weights, goes back to Vojta \cite{Voj89}. It is stated in \cite[\S2.4]{EF13}, where
it follows from a bound for the height of $T(L,c)$ \cite[Prop.~17.5]{EF13}, and in Evertse's
survey \cite[Thm.~C]{Ev10}. In a note on the Harder--Narasimhan filtration of a vector space
with several valuations \cite[\S7]{Ev23}, Evertse builds such a collection from subsets of
adapted bases; its size is of the order of our $J$. Lemma~\ref{lem:types} is an explicit form
of this finiteness.

For the intervals we found no uniform statement. When many weights occur, the literature
discretises them: Evertse and Schlickewei cover the weights by a net and apply the theorem
once per point of the net \cite[\S21]{ES02}, which gives the factor $(9n^2\eps^{-1})^{ns}$ in
\cite[Cor.~3.2]{EF13}, and the same method is used in \cite{NQ25}. Evertse remarks in
\cite[\S6]{Ev10} that the Faltings--W\"ustholz auxiliary polynomial works for solutions
$x_1,\dots,x_M$ ``not necessarily belonging to the same class''. This is the mechanism behind
Theorem~\ref{thm:F}, but it is not stated there as a theorem about intervals. We checked the
works citing \cite{EF13} and \cite{BE08}.

For the complexity bounds, Bugeaud writes that the results of \cite{EF13} give \eqref{eq:BE}
``for $\eta$ in a slightly larger interval than $[0,1/11)$'' \cite[\S12]{Bug15}, without a
value or a proof. We know of no exponent beyond $1/11$ in \eqref{eq:BE}, and no explicit
$\delta$ in the continued fraction analogue \cite[(6.1)]{Bug13}, other than the $1/4$ of
\cite{Ste26b}.

\subsection*{Formalization}
Theorems~A--C are proved in Lean~4 with Mathlib, using only the three standard axioms of Lean.
So are Theorem~\ref{thm:F}, its consequences in Section~\ref{sec:systems}, and everything they use.
The code is in the repository \url{https://github.com/rwst/Subspace-Theorems}, commit
\texttt{30cdebc}. Appendix~\ref{app:lean} lists, for every numbered statement, the Lean
declaration that proves it.

\section{The interval result for a set of weights}\label{sec:F}

We use the notation of \cite[\S2]{EF13}. Let $K$ be a number field with set of places $M_K$,
normalised absolute values $\|\cdot\|_v$, $n\ge2$, and $L=(L_i^{(v)}:v\in M_K,1\le i\le n)$ a
tuple of linear forms with coefficients in $K$ such that
\begin{itemize}
\item $L_1^{(v)},\dots,L_n^{(v)}$ are linearly independent for every $v$;
\item the set $\{L_1,\dots,L_r\}=\bigcup_v\{L_1^{(v)},\dots,L_n^{(v)}\}$ is finite, with
  $r\le R$.
\end{itemize}
Put $\Delta_L=\prod_v\|\det(L_1^{(v)},\dots,L_n^{(v)})\|_v$ and
$H_L=\prod_v\max_{i_1<\dots<i_n}\|\det(L_{i_1},\dots,L_{i_n})\|_v$ \cite[(2.11), (2.12)]{EF13}.
A tuple of weights $c=(c_{iv})$ is \emph{normalised} if $c_{iv}=0$ for all but finitely many
$v$ and
\begin{equation}\label{eq:norm}
  \sum_{i=1}^nc_{iv}=0\quad(v\in M_K),\qquad\sum_{v\in M_K}\max_ic_{iv}\le1,
\end{equation}
which are (2.7)--(2.9) of \cite{EF13}. For a subspace $U\subseteq\Qbar^n$ of dimension $k>0$ the
weight is
\begin{equation}\label{eq:weight}
  w_{L,c}(U)=\sum_{v\in M_K}w_v(U),\qquad
  w_v(U)=\min\Bigl\{\sum_{i\in I}c_{iv}:\#I=k,\ (L_i^{(v)}|_U)_{i\in I}\text{ linearly
  independent}\Bigr\},
\end{equation}
and $w_{L,c}(\{0\})=0$ \cite[(2.19), (2.20)]{EF13}. There is a unique proper subspace
$T=T(L,c)$ such that $w(T)/(n-\dim T)\ge w(U)/(n-\dim U)$ for every proper subspace $U$, with
$\dim T$ minimal subject to this \cite[(2.21) and Lemma~15.2]{EF13}.

\begin{theorem}\label{thm:F}
Let $L$ be as above, $0<\delta\le1$, and $\cC$ a set of normalised tuples of weights. Let
$S_1\subseteq M_K$ be a finite set with $\#S_1\le R^n$, and $\gamma=(\gamma_v)$ reals with
\[
  \gamma_v\ge0,\qquad\sum_v\gamma_v\le1,\qquad\gamma_v=0\text{ for }v\notin S_1,
\]
such that every $c\in\cC$ satisfies $c_{iv}\le\gamma_v$ for all $i,v$ and $c_{iv}=0$ for
$v\notin S_1$. Put $C_0=\max(H_L^{1/R},n^{1/\delta})$ and $\omega_0=\delta^{-1}\log3R$. Then
there are
\begin{itemize}
\item proper subspaces $T_1,\dots,T_J$ of $\Qbar^n$ with $J\le(n+1)\,2^{2^n\#S_1}$, and
\item reals $Q_1,\dots,Q_M\ge C_0$ with $M\le J\,n^{n\#S_1}\,m(n,R,\delta)$,
\end{itemize}
such that for every $c\in\cC$ the subspace $T(L,c)$ is one of the $T_j$, and every real $Q\ge
C_0$ for which
\[
  \{x\in\Qbar^n:H_{L,c,Q}(x)\le\Delta_L^{1/n}Q^{-\delta}\}\not\subseteq T(L,c)
\]
lies in $[Q_h,Q_h^{\omega_0})$ for some $h$. Here $m(n,R,\delta)$ is an explicit integer with
\[
  m(n,R,\delta)\le60n(n+3)\bigl(64000\,n^84^n(1+\log64nR)+2\bigr)\,\delta^{-2}(1+\log\delta^{-1}).
\]
\end{theorem}

For $\cC=\{c\}$ with $S_1$ of size at most $R^n$ this is Theorem~2.3 of \cite{EF13}, with a
different constant in place of $m_0$. Some remarks on the hypotheses follow.
\begin{itemize}
\item The weights of $\cC$ need not be close to each other. They only have to share the bound
  $\gamma$ and the support $S_1$.
\item The condition $\#S_1\le R^n$ is not a restriction. Replacing $R$ by a larger integer only
  changes $m(n,R,\delta)$ through $\log R$. In the applications $\#S_1$ is the number of places
  above $\infty$ and the given primes, a small constant.
\item The bound $J$ depends neither on the forms nor on $\cC$. That the $T(L,c)$ come from a
  finite collection depending only on $L$ is known \cite{Voj89}, \cite[\S2.4]{EF13},
  \cite[\S7]{Ev23}; see the introduction.
\item The intervals are not effectively computable, just as in \cite{EF13}: their number is
  explicit, their location is not.
\end{itemize}

\section{Proof of Theorem~\ref{thm:F}}\label{sec:proof}

The proof follows \cite{EF13}. Theorem~2.3 there is deduced from a semistable case,
Theorem~8.1, in \S\S15--18, and Theorem~8.1 is proved in \S\S9--14. We go through the places
where a single tuple $c$ is used. There are three. The Roth machinery of \S\S13--14 needs a
common bound for the weights (\S\ref{ss:roth}). The exceptional subspaces must be counted
(\S\ref{ss:types}). The reduction of \S18 must be made with one quotient system for many weights
(\S\ref{ss:reduction}). Everything else in \cite{EF13} is applied to one $c$ at a time.

\subsection{The semistable case}\label{ss:roth}
Theorem~8.1 of \cite{EF13} assumes, besides \eqref{eq:norm}, that there is a non-archimedean
place $v_0$ with $c_{i,v_0}=0$ and $L_i^{(v_0)}=X_i$ (8.8), and that $(L,c)$ is semistable,
$T(L,c)=\{0\}$ (8.9). It concludes that the $Q\ge C_2$ with a nonzero point of twisted height at
most $Q^{-\delta}$ lie in $m_2\asymp\delta^{-2}\log\delta^{-1}$ intervals.

\begin{proposition}\label{prop:81}
Let $\cC$ be a set of tuples of weights such that every pair $(L,c)$, $c\in\cC$, satisfies
\textup{(8.1)--(8.9)} of \cite{EF13} with the same place $v_0$. Let $\gamma$ and $S_1$ be as in
Theorem~\ref{thm:F}, with $\gamma_{v_0}=0$. Then the $Q\ge C_2$ for which some $c\in\cC$ has a
nonzero point $x$ with $H_{L,c,Q}(x)\le Q^{-\delta}$ lie in at most $(n-1)(m-1)$ intervals
$[Q',Q'^{2\omega}]$. Here $m\asymp\delta^{-2}\log(R/\delta)$ is the number of blocks of the
auxiliary polynomial and $\omega=4m/\eps$, as in \textup{(14.1)} and \textup{(14.4)} of
\cite{EF13} up to the constants.
\end{proposition}

\begin{proof}[Sketch]
Suppose there are more intervals. Lemma~9.4 of \cite{EF13} is a pigeonhole argument over the
exceptional $Q$; applied to the union of the exceptional sets of all $c\in\cC$, it gives $m$
reals $Q_1<\dots<Q_m$ with $\log Q_{h+1}\ge\omega\log Q_h$, each exceptional for some
$c_h\in\cC$, with a common index $k\in\{1,\dots,n-1\}$ in (14.5). For each $h$, \S11 of
\cite{EF13} gives the points $\hat h_{h1},\dots,\hat h_{h,N-1}$ of the $k$-th exterior power,
with exponents $\hat c^{(h)}$ computed from $c_h$; this is done one block at a time. The rest of
the argument changes as follows.
\begin{itemize}
\item \emph{The auxiliary polynomial} (Lemma~13.5 and Prop.~13.6 of \cite{EF13}). The
  polynomial $P$ has $m$ blocks of variables, and the vanishing conditions in block $h$ are
  read with the exponents $\hat c^{(h)}$. Hoeffding's inequality, Lemma~13.2 of \cite{EF13}, is
  stated for independent random variables with different distributions, and Lemma~13.3 already
  allows one tuple per block. The condition at a place $v\in S_1$ is taken relative to the
  common bound $\gamma_v$ instead of $\max_ic_{iv}$. Evertse and Ferretti replace $S_1$ by a set
  $S_0$ of representatives. Their choice groups the places by the cell of $c_v/\gamma_v$, which
  depends on $c$. We keep all of $S_1$ instead. The counting step then needs
  $2(\#S_1+1)\le e^{m\eps^2/2}$, which holds for the $m$ of (14.1) since $\#S_1\le R^n$.
\item \emph{The final estimate} (\S14 of \cite{EF13}). The upper bound for $|P(\dots)|_v$ at the
  points uses $\gamma_v$ in place of $\max_ic_{iv}$. The normalisation $\sum_v\gamma_v\le1$
  takes the place of (8.4), and the remaining inequalities do not compare different blocks.
\item \emph{The non-vanishing} (Prop.~12.1 of \cite{EF13}) does not involve the weights.
\end{itemize}
The product formula then gives the same contradiction as in \cite[\S14]{EF13}.
\end{proof}

\subsection{Weight types}\label{ss:types}
Let $S_1$ be a finite set of places. The \emph{weight type} of a subspace $U\subseteq\Qbar^n$ is
\[
  \tau(U)=\bigl(\dim U,\ (\mathcal B_v(U))_{v\in S_1}\bigr),\qquad
  \mathcal B_v(U)=\{I\subseteq\{1,\dots,n\}:(L_i^{(v)}|_U)_{i\in I}\text{ is a basis of }U^*\}.
\]
There are at most $(n+1)\,2^{2^n\#S_1}$ weight types.

\begin{lemma}\label{lem:types}
Let $c_1,c_2$ be normalised tuples of weights vanishing outside $S_1$. If $T(L,c_1)$ and
$T(L,c_2)$ have the same weight type, then $T(L,c_1)=T(L,c_2)$. Consequently, for every set $\cC$
of normalised weights vanishing outside $S_1$, the subspaces $T(L,c)$, $c\in\cC$, are at most
$(n+1)\,2^{2^n\#S_1}$ subspaces.
\end{lemma}

\begin{proof}
By \eqref{eq:weight}, $w_v(U)$ is the minimum of $\sum_{i\in I}c_{iv}$ over
$I\in\mathcal B_v(U)$ for $v\in S_1$, and $w_v(U)=0$ for $v\notin S_1$. So $w_{L,c}(U)$ and
$\dim U$ depend only on $\tau(U)$ and $c$. Let $T_j=T(L,c_j)$. For the weights $c_1$, the
subspace $T_2$ has the same dimension and the same weight as $T_1$. So $T_2$ satisfies the
defining condition of $T(L,c_1)$ in \cite[(2.21)]{EF13}, with the same minimal dimension, and
the uniqueness of $T(L,c_1)$ \cite[Lemma~15.2]{EF13} gives $T_2=T_1$.
\end{proof}

\subsection{The reduction}\label{ss:reduction}
Let $c$ be normalised and $T=T(L,c)$ of dimension $k$. In \S18 of \cite{EF13} the general case
is reduced to the semistable one. At each place $v$ one selects index sets $I_v$ of size $n-k$,
by the ordering (18.3) of the $c_{iv}$, and builds from them a tuple of forms $L''$ on
$\Qbar^n/T$ and weights $d$ for which Theorem~8.1 applies. The weights $d$ are bounded by
$\max_ic_{iv}$ at every place and vanish where $c$ vanishes. The forms $L''$ depend on $c$ only
through the selected index sets $I_v$, and on places outside $S_1$, where $c=0$, the selection
does not depend on $c$ at all. The thresholds of \S18 depend on $L$, $n$, $R$ and $\delta$
only, because they come from the height bound for $T$ \cite[Prop.~17.5]{EF13}.

\begin{proof}[Proof of Theorem~\ref{thm:F}]
Split $\cC$ into groups according to the pair $\bigl(T(L,c),(I_v)_{v\in S_1}\bigr)$. By
Lemma~\ref{lem:types} there are at most $J$ subspaces, and for each of them at most $n^{n\#S_1}$
selections. Within one group the quotient forms $L''$ are the same for every member, and the
weights $d$ satisfy the hypotheses of Proposition~\ref{prop:81} with the common bound $\gamma$.
Each member is semistable by the argument of \cite[\S18]{EF13} for that member. So
Proposition~\ref{prop:81} gives one family of intervals for the whole group. Cutting them into
intervals of ratio $\omega_0$, as at the end of \cite[\S18]{EF13}, gives at most
$m(n,R,\delta)$ intervals per group. In all there are at most $J\,n^{n\#S_1}m(n,R,\delta)$.
\end{proof}

\section{Systems of inequalities}\label{sec:systems}

\subsection{Approximation domains}
As in \cite[\S5]{EF13}, a system of inequalities is turned into twisted heights by centring and
normalising its exponents. For a set of systems this has to be done with one normalisation for
all of them.

Let $K$ be a number field, $S_\infty$ its set of archimedean places, $S_{\rm fin}$ a finite set
of finite places, and $t=\#S_\infty+\#S_{\rm fin}$. Let $L_i^{(v)}$ ($v\in S_\infty\cup
S_{\rm fin}$, $1\le i\le n$) be linear forms, linearly independent at each $v$, with at most $s$
distinct forms. For a tuple of exponents $c=(c_{iv})$ and $Q>1$, the \emph{approximation domain}
is
\[
  \Pi(c,Q)=\{x\in K^n:\|L_i^{(v)}(x)\|_v\le Q^{c_{iv}}\ (v\in S_\infty\cup S_{\rm fin}),\
  \|x\|_v\le1\ (v\notin S_\infty\cup S_{\rm fin})\}.
\]
The weight of $c$ is $\sum_v\sum_ic_{iv}$ and its absolute weight is $\sum_v\sum_i|c_{iv}|$,
counted with the local degrees at the archimedean places.

\begin{proposition}\label{prop:33}
Let $0<\eps\le4nA$ and let $\cC$ be a set of exponents of weight at most $-\eps$ such that $t$
times the absolute weight of every $c\in\cC$ is at most $A$. Let $R\ge s+n$ with $t\le R^n$, put
$\delta'=\eps/(4nA)$, and let $X_0$ be at least the explicit threshold, which is linear in the
logarithmic height of the forms. Then there are at most $J=(n+1)\,2^{2^nt}$ proper subspaces $W$
of $K^n$ and at most $J\,n^{nt}\,m(n,R,\delta')$ reals $s\ge X_0$ such that for every $c\in\cC$
and every $Q$ with $\log Q\ge X_0$,
\[
  \Pi(c,Q)\subseteq W\ \text{for one of the }W,\qquad\text{or}\qquad
  s\le\log Q<\delta'^{-1}\log(3R)\cdot s\ \text{for one of the }s.
\]
\end{proposition}

\begin{proof}
For one $c$ this is the deduction of \cite[\S5]{EF13}: the exponents are centred at each place
and divided by a normalising factor $\kappa$, and Theorem~2.3 of \cite{EF13} is applied at the
level $Q^\kappa$. For the set $\cC$ we take one factor $\kappa=2A$. Since $t$ times the absolute
weight is at most $A$, the centred exponents of every $c\in\cC$ are at most $\gamma_v=1/t$ at
each of the $t$ places, and $\sum_v\gamma_v=1$. Theorem~\ref{thm:F} applies with
$S_1=S_\infty\cup S_{\rm fin}$ and the exponent $\delta'$.
\end{proof}

A system with forms over a Galois extension $E$ of $K$ and solutions in $K^n$ is reduced to
Proposition~\ref{prop:33} over $E$ by taking conjugates of the forms, as in \cite[\S5]{EF13}. In
the Lean files this is done for normalised systems with common forms, constants and
normalisation, with $t$ the number of places of $E$ above $\infty$ and above $S$.

\subsection{Real forms}
The applications need forms in $1,\alpha,\alpha^2$ at $\infty$ and the coordinates at finitely
many primes. For $x\in\Z^n$ write $|x|=\max_i|x_i|$, and $|\cdot|_l$ for the $l$-adic absolute
value normalised by $|l|_l=l^{-1}$. The following is the form of the quantitative Subspace
Theorem used in \cite{Ste26b}, now for a set of exponents and in interval form.

\begin{theorem}\label{thm:real}
Let $\alpha$ be a real algebraic number, $S$ a finite set of primes, and $2\le n\le4$. There is
$K=K(\alpha,S,n)\ge1$ with the following property. Let
\begin{itemize}
\item $L_1,\dots,L_n$ be linear forms in $X_1,\dots,X_n$ with coefficients
  $c_0+c_1\alpha+c_2\alpha^2$, where $c_j\in\Z$ and $|c_j|\le B$ for some $B\ge1$, and with
  $|\det(L_1,\dots,L_n)|\ge\gamma$ for some $0<\gamma\le1$;
\item $C\ge1$ and $0<\delta\le1$ be reals;
\item $\cC$ be a set of exponent tuples $(e,f)=((e_i)_i,(f_{l,i})_{l\in S,i})$, each with
  $e_i\le1$, some $e_i=1$, some $e_i\le1-\delta/(2n)$, all $f_{l,i}\le0$, for each $l$ some
  $f_{l,i}=0$, and $\sum_ie_i+\sum_{l\in S}\sum_if_{l,i}\le-\delta$.
\end{itemize}
Then there are at most $K$ proper subspaces of $\Q^n$ and at most $K\delta^{-2}(1+\log\delta^{-1})$
reals $s$ with the following property. Let $(e,f)\in\cC$, and let $x\in\Z^n\setminus\{0\}$ satisfy
\begin{gather*}
  |L_i(x)|\le C|x|^{e_i},\qquad |x_i|_l\le|x|^{f_{l,i}}\qquad(1\le i\le n,\ l\in S),\\
  \log|x|\ge K\delta^{-1}(1+\log B+\log C+\log\gamma^{-1}).
\end{gather*}
Then $x$ lies in one of the subspaces, or $s\le\log|x|<K\delta^{-1}s$ for one of the reals $s$.
\end{theorem}

\begin{proof}[Sketch]
As for one tuple, the forms are rescaled and the system is passed to the Galois closure $E$ of
$\Q(\alpha)$, where it becomes a normalised system over $E$. The rescaled forms
depend on $e$ only through the set $\{i:e_i>1-\delta/(2n)\}$. So $\cC$ falls into at most $2^n$
classes with common forms, constants and normalisation. In each class
Proposition~\ref{prop:33} over $E$ gives at most $J$ subspaces and
$J\,n^{nt}m\ll\delta^{-2}(1+\log\delta^{-1})$ intervals, where $t$ is the number of places of
$E$ above $\infty$ and above $S$. All constants depend only on $\alpha$, $S$ and $n$.
\end{proof}

The constant $K$ is explicit; it is at least $2^nJ$ with $J=(n+1)2^{2^nt}$. For the steering word
($\alpha=0$, $S=\{2,3\}$, $n=3$, $t=3$) this allows up to $2^{29}$ planes. The number is large
but independent of $\delta$, and that is all the applications need. For $\cC=\{(e,f)\}$,
Theorem~\ref{thm:real} is the interval form of the theorem used in \cite{Ste26b}, with one
subspace.

\subsection{Solutions of geometrically growing height}
In interval form, a sequence of solutions whose heights grow fast is controlled without any gap
principle.

\begin{lemma}\label{lem:window}
Let $g_1,g_2,\dots$ be positive reals with $qg_j\le g_{j+1}$ for some $q>1$, let $\rho\ge1$, and
let $\mathcal T$ be a finite set of reals. If every $j$ in a finite set $I$ has $s\le g_j<\rho s$
for some $s\in\mathcal T$, then $\#I\le\#\mathcal T\,(\log\rho/\log q+1)$.
\end{lemma}

\begin{proof}
An interval $[s,\rho s)$ contains at most $\lfloor\log\rho/\log q\rfloor+1$ of the $g_j$.
\end{proof}

\begin{corollary}\label{cor:doubling}
In the situation of Theorem~\ref{thm:real}, let $x_1,x_2,\dots$ be points with $q\log|x_j|\le
\log|x_{j+1}|$ for some $q>1$, each a solution, above the threshold, of the system of some
$(e,f)\in\cC$. Then all but at most
\[
  K\delta^{-2}(1+\log\delta^{-1})\Bigl(\frac{\log(K\delta^{-1})}{\log q}+1\Bigr)
\]
of the $x_j$ lie in the union of at most $K$ proper subspaces of $\Q^n$, which do not depend on
the sequence.
\end{corollary}

\begin{proof}
Theorem~\ref{thm:real} and Lemma~\ref{lem:window} with $g_j=\log|x_j|$ and $\rho=K\delta^{-1}$.
\end{proof}

In the applications the points on the subspaces are counted separately, by a Liouville argument
or by one more application of the Subspace Theorem in fewer variables. Each subspace holds
$O(\log\delta^{-1})$ points.

\section{Integer bases: proof of Theorem~A}\label{sec:A}

Let $b\ge2$ and let $\xi$ be real, algebraic and irrational. We may assume $0<\xi<1$. Suppose
$p(n,\xi,b)\le c\,n(\log n)^u$ for all large $n$, with $0\le u<1/2$, and put $v=(u+1/2)/2$.

\subsection{The points}
By Lemma~6.1 of \cite{BE08} (see \cite[Lemma~3.1]{Ste26b}) there are $C\in\N$, integers
$0\le r_n<t_n$ and $p_n\in\Z$ such that for every $n$
\begin{gather*}
  2t_n<t_{n+1},\qquad t_n\le(2n+2)^{C(n+1)},\qquad |p_n|\le b^{t_n},\qquad
  b\nmid p_n\text{ if }r_n\ge1,\\
  |b^{t_n}\xi-b^{r_n}\xi-p_n|\le(b^{t_n})^{-(\log t_n)^{-v}}.
\end{gather*}
Put $x_n=(b^{t_n},b^{r_n},p_n)$ and $s_n=t_n-r_n$, the \emph{gap}. Let
$\eps=(\log t_{2N})^{-v}$, so that $\eps^{-1}=N^{v+o(1)}$, and count the $N$ points $x_n$,
$N\le n<2N$. A point has a \emph{large gap} if $\eps t_n<4s_n$ and a \emph{small gap}
otherwise.

\subsection{Large gaps}

\begin{proposition}\label{prop:planesA}
There is $K\ge1$ such that for $0<\eps\le1$ there are at most $K$ proper subspaces of $\Q^3$ and
at most $K\eps^{-2}(1+\log\eps^{-1})$ reals $s$ with the following property. Every point
$(b^t,b^r,p)$ with $r<t$, $|p|\le b^t$, $|b^t\xi-b^r\xi-p|\le(b^t)^{-\eps}$ and $K\eps^{-1}\le
t\log b$ lies in one of the subspaces, or satisfies $s\le t\log b<K\eps^{-1}s$ for one of the
reals.
\end{proposition}

\begin{proof}
Let $k=\lfloor2/\eps\rfloor+1$ and put the point into the class $\ell=\lfloor kr/t\rfloor$. For
$l\mid b$ put $g_l=-v_l(b)\log l/\log b$, so that $\sum_{l\mid b}g_l=-1$. As in
\cite[Prop.~3.3]{Ste26b}, a point of class $\ell$ satisfies the system with the forms $X_1$, $X_2$,
$\xi X_1-\xi X_2-X_3$ at $\infty$, the coordinates at the primes $l\mid b$, and the exponents
\[
  e^{(\ell)}=\bigl(1,(\ell+1)/k,-\eps\bigr),\qquad f^{(\ell)}_l=g_l\,(1,\ell/k,0).
\]
These sum to $1/k-\eps=-\delta$ with $\delta\ge\eps/2$. Apply Theorem~\ref{thm:real} once, with
$n=3$, $S$ the primes dividing $b$, and $\cC=\{(e^{(\ell)},f^{(\ell)}):0\le\ell<k\}$. The forms do
not depend on $\ell$.
\end{proof}

In \cite{Ste26b} each class needed its own application, at a cost of about $2/\eps$ times the
count. Here the number of classes does not enter. Each proper subspace of $\Q^3$ contains at
most $K(1+\log\eps^{-1})$ of the points $x_n$ with a large gap \cite[Lemma~3.5]{Ste26b}. Its
proof uses Liouville's inequality when the gap varies on the plane, and a window of ratio
$4/\eps$ when it is constant.

\subsection{Small gaps}
The points $y_n=(b^{t_n}-b^{r_n},p_n)$ with a small gap satisfy a two-variable system with one
tuple of exponents \cite[Lemma~3.6]{Ste26b}. By Theorem~\ref{thm:real} with $\cC$ a single
tuple, they lie on one line or have $t_n\log b$ in $K\eps^{-2}(1+\log\eps^{-1})$ intervals
$[s,K\eps^{-1}s)$. A line contains at most one of them, because $b$ is prime to $b^s-1$
\cite[Lemma~3.7]{Ste26b}.

\subsection{The count}

\begin{proposition}\label{prop:countA}
There is $A>0$ such that for $0<\eps\le1$ and every finite set $I$ of indices whose points satisfy
the conditions of Lemma~6.1 of \cite{BE08} with exponent $\eps$ and $A\eps^{-1}\le t_n\log b$,
\[
  \#I\le A\eps^{-2}(1+\log\eps^{-1})^2.
\]
\end{proposition}

\begin{proof}
The exponents $t_n\log b$ more than double from one index to the next. By
Lemma~\ref{lem:window} with $q=2$ and $\rho=K\eps^{-1}$, the intervals of
Proposition~\ref{prop:planesA} and of the small gaps hold $O(\eps^{-2}(1+\log\eps^{-1})^2)$
points. The at most $K$ subspaces hold $K^2(1+\log\eps^{-1})$ points with a large gap, and the
one line holds at most one point with a small gap.
\end{proof}

\begin{proof}[Proof of Theorem~A]
For large $N$ the $N$ indices $N\le n<2N$ satisfy the hypotheses of
Proposition~\ref{prop:countA}, so $N\le A\eps^{-2}(1+\log\eps^{-1})^2=N^{2v+o(1)}$. This is
false for large $N$ since $2v<1$.
\end{proof}

\section{Continued fractions: proof of Theorem~B}\label{sec:B}

Let $\alpha=[0;a_1,a_2,\dots]$ be algebraic of degree at least three, and suppose
$p(\ell,a)\le C\ell(\log\ell)^\delta$ for all large $\ell$, with $0\le\delta<1/2$. Then
$a_\ell\le M$ for some integer $M\ge1$. We use the points of \cite[\S4]{Ste26b}, which come from
Bugeaud's quadratic approximations \cite[\S3]{Bug13}. A repetition in the word $a$ gives a point
$v=(v_0,\dots,v_3)\in\Z^4$, the coefficient vector of the minimal polynomial of a quadratic number,
of height $H=q_wq_{w+r}$, together with $Q=q_w$. For the forms
\[
  L_1=\alpha^2X_0-\alpha(X_1+X_2)+X_3,\qquad L_2=\alpha X_0-X_1,\qquad L_3=\alpha X_0-X_2
\]
these points satisfy the inequalities of \cite[Lemma~4.1]{Ste26b}. Put $\rho=\log Q/\log H\in
[0,1/2]$ and let $T_0=\{X_1=X_2\}$.

\begin{proposition}\label{prop:fourB}
There is $K\ge1$ such that for $0<\eps\le1$ there are at most $K$ proper subspaces of $\Q^4$ and
at most $K\eps^{-2}(1+\log\eps^{-1})$ reals $s$ with the following property. Every point of
\cite[Lemma~4.1]{Ste26b} with exponent $\eps$ and $\log H\ge K\eps^{-K}$ lies in one of the
subspaces or has $s\le\log H<K\eps^{-1}s$ for one of the reals.
\end{proposition}

\begin{proof}
Cover $[0,1/2]$ by cells $[\rho^-,\rho^+]$ of length at most $\eps/4$. On one cell the forms
$L_1,L_2,L_3,X_0$, of determinant $-1$, are bounded at $v$ by $24(M+1)^4|v|^{e_i}$ with
\[
  e=(-1-\eps,\ 2\rho^+-1,\ 1-2\rho^-,\ 1),\qquad\textstyle\sum_ie_i\le-\eps/2
\]
\cite[Prop.~4.3]{Ste26b}. Apply Theorem~\ref{thm:real} once, with $n=4$, $S=\emptyset$,
$\delta=\eps/2$ and $\cC$ the set of the exponents of all cells.
\end{proof}

\begin{proposition}\label{prop:countB}
There is $K$ with the following property. Let $0<\eps\le1$ and let $v_k$ $(k\in I)$ be points as
above with heights $H_k$ such that $\log H_k\ge1$, $\log H_{k+1}\ge2\log H_k$ and
$\log H_k\ge K\eps^{-K}$. Then $\#I\le K\eps^{-2}(1+\log\eps^{-1})^2$.
\end{proposition}

\begin{proof}
By Proposition~\ref{prop:fourB} and Lemma~\ref{lem:window} with $q=2$, the points with $\log H$
in one of the intervals number $O(\eps^{-2}(1+\log\eps^{-1})^2)$. The remaining points lie in
at most $K$ proper subspaces. Each proper subspace contains $O(\log\eps^{-1})$ points outside
$T_0$, by a three-variable application at the fixed exponent $3/4$ and a window
\cite[Prop.~4.3, steps (ii) and (iii)]{Ste26b}. On $T_0$ the forms $L_1,L_2,X_0$ in
$(X_0,X_1,X_3)$ have determinant $1$ and exponents $(-1-\eps,0,1)$. Theorem~\ref{thm:real} with a
single tuple puts the points of $T_0$ into one plane or into $K\eps^{-2}(1+\log\eps^{-1})$
intervals. A plane holds $O(\log\eps^{-1})$ points \cite[Prop.~4.3, step (ii)]{Ste26b}, and
Lemma~\ref{lem:window} bounds the points in the intervals.
\end{proof}

\begin{proof}[Proof of Theorem~B]
This is the proof of \cite[Thm.~1.2]{Ste26b} with Proposition~\ref{prop:countB} in place of
\cite[Prop.~4.3]{Ste26b}. Take lengths $\ell_k$ with $\log\ell_k\asymp k\log k$. The
configurations of \cite[Lemma~4.2]{Ste26b} give points $v_k$ of doubly exponentially growing
height. With $\eps_N^{-1}\asymp(N\log N)^\delta$, Proposition~\ref{prop:countB} gives
\[
  N\le K\eps_N^{-2}(1+\log\eps_N^{-1})^2=(N\log N)^{2\delta+o(1)},
\]
which fails for large $N$ since $2\delta<1$.
\end{proof}

\section{The steering word: proof of Theorem~C}\label{sec:C}

Write $(3/2)^n=m_n+\eps_n$ with $m_n$ the nearest integer and $t_n=2m_{n+1}-3m_n$. A repetition
of length $L$ at $a<c$ forces $|\eps_c-\eps_a|\le(2/3)^L$ \cite[Lemma~2.7]{St26}. Put
$\lambda=\log2/\log3$. Suppose $\pT(k)\le c\,k(\log k)^u$ for all large $k$, with $0\le u<1/2$,
and put $v=(u+1/2)/2$. As in \cite[\S5]{Ste26b} there are a doubling sequence $2\tau_n<\tau_{n+1}$
and pairs $(a_n,c_n)=(r_n+2,\tau_n+2)$ with $|\eps_{c_n}-\eps_{a_n}|\le(3^{c_n})^{-\eps}$. For
such a pair let $s=c-a$ be its gap, $M=m_c-m_a$, and $x=(3^c,2^s3^a,2^cM)\in\Z^3$ its kernel
point.

\begin{proposition}\label{prop:planesC}
There is $K\ge1$ such that for $0<\eps\le1$ there are at most $K$ proper subspaces of $\Q^3$ and
at most $K\eps^{-2}(1+\log\eps^{-1})$ reals $s$ with the following property. The kernel point of
every pair $2\le a<c$ with $|\eps_c-\eps_a|\le(3^c)^{-\eps}$ and $K\eps^{-1}\le c\log3$ lies in
one of the subspaces, or $s\le c\log3<K\eps^{-1}s$ for one of the reals.
\end{proposition}

\begin{proof}
Let $k=\lfloor2/\eps\rfloor+1$ and sort the pairs by $\ell=\lfloor ks/c\rfloor$. The pairs of
class $\ell$ satisfy the system of \cite[Prop.~5.2]{Ste26b} with the forms $X_0-X_1-X_2,X_1,X_2$ at
$\infty$, the coordinates at $2$ and $3$, and the exponents
\[
\begin{array}{c|ccc}
 & X_0-X_1-X_2 & X_1 & X_2\\\hline
 \infty & \lambda-\eps & 1-(1-\lambda)\ell/k & 1\\
 2 & 0 & -\lambda\ell/k & -\lambda\\
 3 & -1 & -(1-(\ell+1)/k) & 0
\end{array}
\]
These sum to $1/k-\eps=-\delta$ with $\delta\ge\eps/2$. Apply Theorem~\ref{thm:real} once, with
$\alpha=0$, $n=3$, $S=\{2,3\}$ and $\cC$ the set of the $k$ exponent tables.
\end{proof}

\begin{proposition}\label{prop:countC}
There is $A>0$ such that for $0<\eps\le1$, every doubling sequence $2\tau_n<\tau_{n+1}$ with
$r_n<\tau_n$, and every finite set $I$ of indices with $|\eps_{c_n}-\eps_{a_n}|\le
(3^{c_n})^{-\eps}$ and $A\eps^{-1}\le c_n\log3$,
\[
  \#I\le A\eps^{-2}(1+\log\eps^{-1})^2.
\]
\end{proposition}

\begin{proof}
Split by the gap as in \cite[\S5]{Ste26b}.
\begin{itemize}
\item \emph{Large gaps}, $\eps c_n<4s_n$. By Proposition~\ref{prop:planesC} the kernel points lie
  in at most $K$ planes or have $c_n\log3$ in one of $K\eps^{-2}(1+\log\eps^{-1})$ intervals. A
  plane carries at most $12(1+\log\eps^{-1})$ indices with a large gap
  \cite[Lemma~5.6]{Ste26b}.
\item \emph{Small gaps}. The points $y=(3^c-2^s3^a,2^cM)$ satisfy a two-variable system with one
  tuple of exponents \cite[Lemma~5.7]{Ste26b}. By Theorem~\ref{thm:real} they lie on one line,
  which carries at most three indices, or have $c_n\log3$ in $K\eps^{-2}(1+\log\eps^{-1})$
  intervals.
\end{itemize}
The positions $c_n=\tau_n+2$ satisfy $c_{n+1}\ge\frac32c_n$. Lemma~\ref{lem:window} with
$q=3/2$ bounds the indices in the intervals by $O(\eps^{-2}(1+\log\eps^{-1})^2)$.
\end{proof}

\begin{proof}[Proof of Theorem~C]
Let $\eps_0=(\log\tau_{2N})^{-v}$ and $\eps=\eps_0/6$. For $N\le n<2N$ the pairs $(a_n,c_n)$
satisfy the hypotheses of Proposition~\ref{prop:countC} for large $N$, so
$N\le A\eps^{-2}(1+\log\eps^{-1})^2=N^{2v+o(1)}$. This contradicts $2v<1$.
\end{proof}

\section{Remarks}

\begin{remark}\label{rem:third}
Theorem~\ref{thm:real} for one tuple at a time, applied once per class, gives
$\eps^{-1}\cdot\eps^{-2}(1+\log\eps^{-1})^2$ in each of the three counts, and so every exponent
below $1/3$ in Theorems~A--C. This intermediate step is also machine-verified
(Appendix~\ref{app:lean}). The exponents form the sequence
\[
  \frac1{11}=\frac1{1+7+3}\ \text{\cite{BE08}},\qquad\frac14=\frac1{1+3}\ \text{\cite{Ste26b}},
  \qquad\frac13=\frac1{1+2},\qquad\frac12.
\]
In each denominator the $1$ is the power of $\eps^{-1}$ spent on the classes. In \cite{BE08}, $7$
counts the subspaces of \cite{ES02} and $3$ the points per subspace; in \cite{Ste26b}, $3$ counts
the subspaces of \cite[Thm.~3.1]{EF13}; in the interval form, $2$ counts the intervals. With
Theorem~\ref{thm:F} the classes cost nothing.
\end{remark}

\begin{remark}\label{rem:limit}
The exponent $1/2$ is the limit of the method. The intervals of Theorem~\ref{thm:F} come from the
$m\asymp\eps^{-2}\log\eps^{-1}$ blocks of the auxiliary polynomial: Hoeffding's inequality needs
$m\eps^2$ large for the index conditions to cut out a positive proportion of the monomials. This
is the $\delta^{-2}$ of the quantitative Roth theorem \cite[Prop.~1.2]{Ev10}; Evertse
describes the dependence on $\delta$ in \cite[Thm.~3.1]{EF13} as ``almost as good as that in
the quantitative Roth's Theorem'' \cite[\S2]{Ev10}. A sequence of doubly
exponentially growing height can put one point into each interval, so a bound of the shape of
Theorem~\ref{thm:F} never gives fewer than $\delta^{-2}$ points. Going beyond $1/2$ needs
$o(\delta^{-2})$ intervals for lacunary sequences of solutions, or a stronger input from the
complexity side than one point with $\eps=(\log n)^{-v}$ per doubling step.
\end{remark}

\begin{remark}
Theorem~\ref{thm:F} applies whenever solutions are sorted into classes with common forms and
different weights. For example, in the bound of \cite[Cor.~3.2]{EF13} for the product
inequality \cite[(3.11)]{EF13}, the factor $(9n^2\eps^{-1})^{ns}$ comes from a grid of weights
for the same forms \cite[\S21]{ES02}. Theorem~\ref{thm:F} gives an interval result for (3.11)
without this factor; a bound for the number of subspaces would still need a gap principle.
Bugeaud's extension of the Cugiani--Mahler theorem to simultaneous approximation \cite{Bug07}
is another application of the quantitative Subspace Theorem to sequences of solutions. We have
not worked out these applications.
\end{remark}

\section*{Acknowledgements}
We thank Yann Bugeaud for drawing our attention to this question. This work relies on a complete
formalization of the Subspace Theorem and its quantitative versions, which is publicly available
and is being incorporated into the TauCeti library.

\appendix
\section{The Lean formalization}\label{app:lean}

Everything is in the repository
\begin{center}\url{https://github.com/rwst/Subspace-Theorems}, commit \texttt{30cdebc},\end{center}
a Lean~4 project on Mathlib. Every declaration below is proved without \texttt{sorry}. Its axioms,
as reported by \texttt{\#print axioms}, are the three standard axioms \texttt{propext},
\texttt{Classical.choice} and \texttt{Quot.sound}. The repository declares no axioms of its own.
The directories are:
\begin{itemize}
\item \texttt{QuantitativeSubspace/} for Sections~\ref{sec:F}--\ref{sec:systems};
\item \texttt{DiophantineApproximation/} for Lemma~\ref{lem:window} and the reduction of systems;
\item \texttt{BugeaudEvertse2008/Additions/} for Section~\ref{sec:A};
\item \texttt{Bugeaud2013/Additions/} for Section~\ref{sec:B};
\item \texttt{Stephan2026a/} for Section~\ref{sec:C}.
\end{itemize}

In the Lean statements:
\begin{itemize}
\item Theorem~\ref{thm:F} is stated for points in $\Omega^n$ for any field extension $\Omega$
  of $K$, with the twisted height of \cite{EF13}. It asks for $\sum_v\max_ic_{iv}\le1$
  separately, although this follows from the bounds $\gamma_v$.
\item Theorem~A is stated for $\xi=\sum_{j\ge0}a_jb^{-j-1}$ with digits $a_j\in\{0,\dots,b-1\}$,
  irrational and algebraic. Theorem~B counts the factors of $a_1a_2\cdots$ as a cardinal,
  possibly infinite. Theorem~C is stated for the word $t_n=2m_{n+1}-3m_n$ with $m_n$ the rounded
  value of $(3/2)^n$.
\item Corollary~\ref{cor:doubling} is not a separate declaration. Each of the three counts
  combines Theorem~\ref{thm:real} with Lemma~\ref{lem:window} directly.
\end{itemize}

In the table, \lean{FS.}\ stands for the namespace \lean{NumberField.FormSystem.}

\medskip
\begin{center}
\scriptsize
\setlength{\tabcolsep}{4pt}
\begin{tabular}{lll}
\hline
Statement & Lean declaration & File\\\hline
\multicolumn{3}{l}{\emph{Sections~\ref{sec:F}--\ref{sec:systems}: \lean{QuantitativeSubspace/}}}\\
Thm.~\ref{thm:F} & \lean{FS.exists\_intervals\_cover\_family} & \lean{FormIntervalFamily}\\
\quad $m(n,R,\delta)$ & \lean{Real.intervalCount\_le\_of\_inv\_le} & \lean{RealFormInterval}\\
Prop.~\ref{prop:81} & \lean{FS.exists\_intervals\_family} & \lean{FormSemistableGap}\\
\quad Lemma 13.5 & \lean{FS.exists\_auxiliaryPolynomial} & \lean{FormAuxiliaryPolynomial}\\
\quad Prop.~13.6 & \lean{FS.exists\_auxiliaryPolynomial\_hasseDeriv\_of\_places} & \lean{FormAuxiliaryPolynomial}\\
\quad \S14 & \lean{FS.false\_of\_chain\_family} & \lean{FormIntervalTheorem}\\
Lemma~\ref{lem:types} & \lean{FS.eq\_of\_isDestabilizing\_of\_weightType\_eq} & \lean{FormWeightType}\\
 & \lean{FS.exists\_finset\_isDestabilizing} & \lean{FormWeightType}\\
\S\ref{ss:reduction} & \lean{FS.quotSystem\_congr} & \lean{FormIntervalFamily}\\
 & \lean{FS.exists\_intervals\_cover\_of\_isDestabilizing\_family} & \lean{FormIntervalFamily}\\
Prop.~\ref{prop:33} & \lean{NumberField.exists\_forall\_mem\_interval\_approxDomain\_ef\_family} & \lean{FormSystemFamily}\\
\quad systems & \lean{NumberField.exists\_forall\_mem\_interval\_of\_efSystemThreshold\_le\_family} & \lean{FormSystemFamily}\\
Thm.~\ref{thm:real} & \lean{Real.subspaceIntervalFamilyAt\_of\_isAlgebraic} & \lean{RealFormFamily}\\
\quad explicit $K$ & \lean{Real.subspaceIntervalFamilyWith\_of\_mem} & \lean{RealFormFamily}\\
\quad $\alpha=0$ & \lean{Real.subspaceIntervalFamilyAt\_zero} & \lean{RealFormFamily}\\
Lemma~\ref{lem:window} & \lean{Real.card\_le\_of\_forall\_mem\_window} & \lean{DoublingSequences}\\
\multicolumn{3}{l}{\emph{Section~\ref{sec:A}: \lean{BugeaudEvertse2008/}}}\\
Thm.~A & \lean{Real.frequently\_lt\_complexity\_of\_lt\_half} & \lean{Additions/Theorem21Half}\\
Prop.~\ref{prop:planesA} & \lean{Real.exists\_plane\_intervals\_family\_of\_approx} & \lean{Additions/Theorem21Half}\\
Prop.~\ref{prop:countA} & \lean{Real.card\_le\_of\_approx\_half} & \lean{Additions/Theorem21Half}\\
\quad planes & \lean{Real.card\_mem\_plane\_le\_log} & \lean{Additions/PlaneLiouville}\\
\quad small gaps & \lean{Real.exists\_line\_intervals\_of\_small\_gap} & \lean{Additions/Theorem21Third}\\
 & \lean{Int.dvd\_of\_pow\_sub\_pow\_mul\_eq} & \lean{Additions/SmallGap}\\
\multicolumn{3}{l}{\emph{Section~\ref{sec:B}: \lean{Bugeaud2013/Additions/}}}\\
Thm.~B & \lean{Nat.frequently\_lt\_encard\_factors\_of\_lt\_half} & \lean{CountHalf}\\
Prop.~\ref{prop:fourB} & \lean{Real.exists\_cover\_four\_intervals\_family} & \lean{CountHalf}\\
Prop.~\ref{prop:countB} & \lean{Real.card\_le\_of\_cfPoints\_half} & \lean{CountHalf}\\
\quad $T_0$ & \lean{Real.exists\_planes\_t0\_intervals} & \lean{CountThird}\\
\quad subspaces & \lean{Real.card\_mem\_le\_of\_ne\_top} & \lean{Subspaces}\\
\quad final step & \lean{Real.cfCount\_one} & \lean{CountHalf}\\
 & \lean{Nat.frequently\_lt\_encard\_factors\_of\_cfCount} & \lean{Delta}\\
\multicolumn{3}{l}{\emph{Section~\ref{sec:C}: \lean{Stephan2026a/}}}\\
Thm.~C & \lean{SteeringWord.frequently\_lt\_complexity\_half} & \lean{QuantitativeHalf}\\
Prop.~\ref{prop:planesC} & \lean{SteeringWord.exists\_plane\_intervals\_family\_of\_violators} & \lean{QuantitativeHalf}\\
Prop.~\ref{prop:countC} & \lean{SteeringWord.card\_le\_of\_violators\_half} & \lean{QuantitativeHalf}\\
\quad planes & \lean{SteeringWord.card\_mem\_plane\_le} & \lean{QuantitativePlaneStep}\\
\quad small gaps & \lean{SteeringWord.exists\_line\_intervals\_of\_violators} & \lean{QuantitativeThird}\\
\multicolumn{3}{l}{\emph{Remark~\ref{rem:third}}}\\
Thm.~A, $1/3$ & \lean{Real.frequently\_lt\_complexity\_of\_lt\_third} & \lean{Additions/Theorem21Third}\\
Thm.~B, $1/3$ & \lean{Real.cfCount\_two} & \lean{CountThird}\\
Thm.~C, $1/3$ & \lean{SteeringWord.frequently\_lt\_complexity\_third} & \lean{QuantitativeThird}\\
\hline
\end{tabular}
\end{center}

\end{document}